\documentclass[11pt, a4paper,twoside]{article}
\usepackage{hyperref}
\usepackage{enumerate}
\usepackage{amsmath}

\usepackage{pb-diagram}
\usepackage[english]{babel}
\usepackage{amssymb,latexsym}
\def \bui#1#2{\mathrel{\mathop{\kern 0pt#1}\limits^{#2}}}

\let\<\langle
\let\>\rangle
\newcommand{\R}{{\mathbb R}}

\newcommand{\lquot}[2]{\raisebox{-0.5ex}{$#2$}\backslash\!\raisebox{0.5ex}{$#1$}}

\newtheorem{example}{Examples}[section]
\newtheorem{thm}{Theorem}[section]
\newtheorem{lemma}[thm]{Lemma}

\newtheorem{cor}[thm]{Corollary}
\newtheorem{remark}[thm]{Remark}
\newtheorem{remarks}[thm]{Remarks}
\newtheorem{definition}[thm]{Definition}
\newtheorem{notation}[thm]{Notation}
\newtheorem{exabout:ample}[thm]{Example}

\begin{document}
\title{Eigenvalue Estimate for the basic Laplacian on manifolds with foliated boundary, part II}
\author{Fida El Chami\footnote{Lebanese University, Faculty of Sciences II, Department of Mathematics, P.O. Box 90656 Fanar-Matn, Lebanon,
		E-mail: \texttt{fchami@ul.edu.lb}},\, Georges Habib\footnote{ Laboratoire de Math\'ematiques et Mod\'elisation, PR2N / EDST,		
		Lebanese University, Faculty of Sciences II, P.O. Box 90656 Fanar-Matn, Lebanon, 
		E-mail: \texttt{ghabib@ul.edu.lb}},\, Ola Makhoul \footnote{ Lebanese University, Faculty of Sciences II, Department of Mathematics, P.O. Box 90656 Fanar-Matn, Lebanon, E-mail: \texttt{ola.makhoul@ul.edu.lb}}, \, Roger Nakad \footnote{Notre Dame University-Louaiz\'e, Faculty of Natural and Applied Sciences, Department of Mathematics and Statistics, P.O. Box 72, Zouk Mikael, Lebanon, E-mail: \texttt{rnakad@ndu.edu.lb}}}

	\date{}





\maketitle 
\begin{abstract}
In \cite{EHGN1}, we gave a sharp lower bound for the first eigenvalue of the basic Laplacian acting on basic $1$-forms defined on a compact manifold whose boundary is endowed with a Riemannian flow. In this paper, we extend this result to the case of 
 basic $p$-forms for $p>1$. As in \cite{EHGN1}, the limiting case allows to characterize the manifold $\lquot{\mathbb{R} \times
B'}{\Gamma}$ for some group $\Gamma$, and where $B'$ denotes the unit closed ball. In particular, we describe the Riemannian product $\mathbb{S}^1\times \mathbb{S}^n$ as the boundary of a manifold.  
\end{abstract}

{\bf Key words}: Riemannian flow, manifolds with boundary, basic Laplacian, eigenvalue, second fundamental form, O'Neill tensor, basic Killing forms, rigidity results.

{\bf Mathematics Subject Classification}: 53C12, 53C24, 58J50, 58J32.

\section{Introduction}

\noindent On a compact manifold $N$ whose boundary $M$ carries a Riemannian flow given by a unit vector field $\xi$ (see Section \ref{sec:2} for the definition), we derive in \cite{EHGN1} a spectral inequality for the basic Laplacian. In fact, by a suitable extension of any basic closed $p$-form to the whole manifold $N$ (see Lemma \ref{lem:1}), we estimated the first eigenvalue of the basic Laplacian (restricted to closed forms) in terms of the principal curvatures when $p=1$ \cite[Thm. 1.1]{EHGN1}. The main tool in our estimate is the so-called {\it Reilly formula} \cite{RS} 
obtained by integrating the Bochner-Weitzenb\"{o}ck formula over $N$ and using the Stokes formula. As a consequence of our estimate, we obtain  several rigidity theorems which mainly characterize the flow as a local product and the boundary as an $\eta$-umbilical submanifold (see \cite[Sec. 5]{EHGN1}). These results can be seen as a foliated version of the work of Raulot and Savo in \cite{RS}.\\

\noindent In this paper, we generalize the results stated in \cite{EHGN1} to any basic $p$-form when $p>1.$ Some of the techniques used in this work are similar to those of $p=1$ (Lemmas \ref{lem:1} and  \ref{lem:2}) but new results are needed to get the estimate and to study the equality case (Lemmas \ref{lem:4}, \ref{lem:geo1}, \ref{lem:geo2} and \ref{lem:geo3}). First, we prove

\begin{thm}\label{thm:main}
Let $(N^{n+2},g)$ be a Riemannian manifold whose boundary $M$ has a positive curvature operator and that $\sigma_{n+1-p}(M)>0$ for some $2\leq p\leq \frac{n}{2}$. Assume that $M$ is endowed with a minimal Riemannian flow given by a unit vector field $\xi.$ Then
\begin{eqnarray}\label{ineq:main}
\lambda'_{1,p}+4 c^2 [\frac{n}{2}]\geq \Big(\sigma_{p+1}(M)-\mathop{{\rm sup}}\limits_Mg(S(\xi),\xi)\Big)\sigma_{n+2-p}(M),
\end{eqnarray}
where $c$ denotes the supremum over $M$ of $|\Omega|$.
\end{thm}

\noindent Inequality \eqref{ineq:main} differs from the one in \cite[Thm. 1.1]{EHGN1} because of the additional term in $c$ that comes from computing the norm of the interior product of the O'Neill tensor with any basic $p$-form (see Lemma \ref{lem:4}). When equality is realized in \eqref{ineq:main}, the term in $c$ turns out to be zero. This yields to the following characterization: 

\begin{thm} \label{thm:mainestimate} Under the same assumptions as in Theorem \ref{thm:main} with $n$ even, if $\sigma_1(M)\geq 0$ and the equality case is realized, the manifold $N$ is then isometric to the quotient $\lquot{\mathbb{R} \times B'}{\Gamma}$, where $B'$ is the unit
closed ball in $\mathbb{R}^{n+1}$ and $\Gamma$ is a co-compact subgroup.
\end{thm}

\noindent As a consequence of this last theorem, we describe manifolds whose boundaries are isometric to the Riemannian product $\mathbb{S}^1 \times \mathbb{S}^n$ and get the following rigidity result:
 
\begin{cor}\label{rigidity1} Let $N$ be an $(n+2)$-dimensional compact
manifold with non-negative curvature operator. Assume that the
boundary $M$ is $\mathbb{S}^1 \times \mathbb{S}^n$ and
$(n+1-p)\sup_M g(S(\xi),\xi) + 4c^2 [\frac n2]\leq 0$. If the inequality $\sigma_{p+1}(M)\geq p$
holds, the manifold $N$ is isometric to $\mathbb{S}^1\times B'$.
\end{cor}

\noindent When changing the sign of the expression in Corollary \ref{rigidity1}, another rigidity result can be obtained:

\begin{cor} \label{rigidity23} Let $N$ be a $(n+2)$-dimensional compact
manifold with non-negative curvature operator. Assume that the
boundary $M$ is $\mathbb{S}^1 \times \mathbb{S}^n$ and
$(n+1-p)\sup_M g(S(\xi),\xi) + 4c^2 [\frac n2]\geq 0$. If the inequality $\sigma_{p+1}(M)\geq p$
holds, the manifold $N$ is isometric to $\mathbb{S}^1\times B'$.
\end{cor}

\noindent Note that we end this paper by an appendix section where we put some technical formulas that were useful in our computations.

\section{Riemannian flows and manifolds with boundary} \label{sec:2}

\noindent Throughout this section, we recall the main ingredients of Riemannian flows defined on a manifold and the basic facts on manifolds with boundary. We will use the preliminaries and the notations of  \cite{EHGN1}.\\

\noindent Let $(M,g)$ be a Riemannian manifold and $\xi$ be a smooth unit vector field on $M$ defining the structure of a Riemannian flow on $M$. In other words, the integral curves of $\xi$ foliate the manifold $M$ such that the leaves are locally equidistant \cite{T}. One can easily see that 
 the endomorphism $h:=\nabla^M\xi$ (known as the O'Neill tensor \cite{O}) defines a skew-symmetric tensor field on the normal bundle $Q=\xi^\perp$. Hence we can associate to $h$ the differential $2$-form $\Omega(Y,Z)=g(h(Y),Z)$ for any sections $Y,Z$ in $\Gamma(Q).$ On the other hand, the normal bundle carries a covariant derivative $\nabla$ compatible with the induced metric $g$ \cite{T}. It is easy to check that the corresponding Levi-Civita connections on $M$ and $Q$ are related for all sections $Z,W$ in $\Gamma(Q)$ via the Gauss-type formulas:
\begin{equation*}
\left\{\begin{array}{ll}
\nabla^M_Z W=\nabla_Z W-g(h(Z),W)\xi, &\textrm {}\\\\
\nabla^M_\xi Z=\nabla_\xi Z+h(Z)-\kappa(Z)\xi,&\textrm {}
\end{array}\right.
\end{equation*}
where $\kappa:=\nabla^M_\xi \xi$ is the mean curvature of the flow. Recall now that a basic form is a differential form $\varphi$ on $M$ such that $\xi\lrcorner \varphi=0$ and $\xi\lrcorner d^M\varphi=0$. We will denote by $\Omega_B(M)$ the set of all such forms. Clearly, basic forms are preserved by the exterior derivative $d^M$ and therefore we can set $d_b:=d^M|_{\Omega_B(M)}.$ For a compact manifold $M$, we consider the $L^2$-adjoint of $d_b$, denoted by $\delta_b$, and define the basic Laplacian as $\Delta_b=d_b\delta_b+\delta_bd_b.$ From the spectral theory of transversally elliptic operators, the basic Laplacian has a discrete spectrum \cite{ElG1,ElG2}.\\

\noindent Since we are going to define Riemannian flows on manifolds with boundary, we then need to recall some basic facts on such manifolds. For this purpose, we let $(N^{n+1},g)$ be a Riemannian manifold of dimension $n+1$ with boundary $M$.  We denote by $\eta_1,\cdots,\eta_{n+1}$ the principal curvatures of $M$ and arrange them so that $\eta_1\leq \eta_2\leq\cdots \leq \eta_{n+1}.$ For any $p\in \{1,\cdots,n+1\},$ we denote the lowest $p$-curvatures $\sigma_p$ by $\sigma_p(x)=\eta_1(x)+\cdots+\eta_{p}(x).$ It is a clear fact that for $p\leq q$, the inequality $\frac{\sigma_p}{p}\leq \frac{\sigma_q}{q}$ holds where the optimality is achieved if and only if either $\eta_1=\eta_2=\cdots=\eta_q$ or $p=q$. Now let $\nu$ be the inward unit nomal vector field on $M$. The shape operator (or the Weingarten tensor) is defined for all $X\in \Gamma(TM)$ as $S(X)=-\nabla^N_X\nu$ where $\nabla^N$ is the Levi-Civita connection of $N$.  Recall the Gauss-Codazzi equation for any $X,Y\in \Gamma(TM),$  
\begin{eqnarray}\label{eq:codazzi}
(\nabla^M_XS)(Y)-(\nabla^M_YS)(X)=R^N(Y,X)\nu.
\end{eqnarray}
As mentionned in \cite{RS}, the Weingarten map admits a canonical extension to any $p$-form $\varphi$ on $M$ by the following: 
$$S^{[p]}(\varphi)(X_1,\cdots,X_p)=\sum_{i=1}^{p}\varphi(X_1,\cdots,S(X_i),\cdots,X_p).$$

The eigenvalues of $S^{[p]}$ are exactly the $p$-curvatures and that means the following inequality 
\begin{equation} \label{eq:1}
\langle S^{[p]}(\varphi),\varphi\rangle \geq \sigma_p(M) |\varphi|^2
\end{equation}
holds, where $\sigma_p(M)$ is the infimum over $M$ of the lowest $p$-curvatures $\sigma_p.$ 

Now, we recall the Reilly formula established in \cite{RS}. For this, we denote by $J^*$ the restriction of differential forms on $N$ to the boundary $M$. At any point $x\in M$, the relation $|J^*\alpha|^2+|\nu\lrcorner \alpha|^2=|\alpha|^2$ is true for any differential form $\alpha$ in $\Lambda^p(N).$ The formula is the following $$
\int_N |d^N\alpha|^2+|\delta^N\alpha|^2=\int_N|\nabla^N\alpha|^2+\langle W_N^{[p]}(\alpha),\alpha\rangle+2\int_M\langle \nu\lrcorner \alpha,\delta^M(J^*\alpha)\rangle\\
+\int_M\mathcal{B}(\alpha,\alpha) $$
where
\begin{eqnarray*}
\mathcal{B}(\alpha,\alpha)&=&\langle S^{[p]}(J^*\alpha),J^*\alpha\rangle+\langle S^{[n+2-p]}(J^*(*_N\alpha)),J^*(*_N\alpha)\rangle\\
&=&\langle S^{[p]}(J^*\alpha),J^*\alpha\rangle+(n+1)H|\nu\lrcorner\alpha|^2-\langle S^{[p-1]}(\nu\lrcorner\alpha),\nu\lrcorner\alpha\rangle,
\end{eqnarray*}
and $W_N^{[p]}$ is the curvature term that appears in the Bochner-Weitzenb\"ock formula for the Laplacian on $N.$ We mention here that $J^*(*_N\alpha)$ is equal (up to a sign) to $*_M(\nu\lrcorner \alpha)$ and that when the curvature operator of $N$ is non-negative, the term $W_N^{[p]}\geq 0.$ 

\noindent Finally, the following boundary problem will be of interest in our study. In fact, given any $p$-form $\varphi$ on $M$, the solution of
\begin{equation} \label{eq:23}
\left\{\begin{array}{ll}
\Delta^N \hat{\varphi}=0 &\textrm{on $N$},\\\\
J^*\hat\varphi=\varphi,\,\, J^*(\delta^N\hat\varphi)=0 &\textrm{on $M$}
\end{array}\right.
\end{equation}
is unique on $N$ by Lemma 3.5.6 in \cite{S}. Moreover, the $p$-form $\hat\varphi$ is co-closed on $N$ and $d^N\hat\varphi\in H^{p+1}(N)$ (see \cite[Lemma 3.1]{BS} for more details).  

\section{Eigenvalue estimate for the basic Laplacian on manifolds with foliated boundary}

In this section, we are going to prove Theorem \ref{thm:main}. 
For this purpose, we need to state the following lemmas already proved in \cite{EHGN1}.

\begin{lemma}\label{lem:1} {\rm  \cite{EHGN1}}
Let $(N^{n+2},g)$ be a Riemannian manifold with boundary $M$ with $W_N^{[p+1]}\geq 0$ for some $1\leq p\leq n$. Assume that $M$ carries a Riemannian flow given by a unit vector field $\xi$. Given a basic closed $p$-form  $\varphi$ and if $\sigma_{n+1-p}(M)>0$ the corresponding solution $\hat\varphi$ of Problem $(\ref{eq:23})$ is then closed and co-closed on $N$.
\end{lemma}

\begin{lemma}\label{lem:2} {\rm \cite{EHGN1}}
Let $(N^{n+2},g)$ be a Riemannian manifold with boundary $M$. Assume that $M$ carries a Riemannian flow given by a unit vector field $\xi$. We have 
$$\langle S^{[p]}(\varphi),\varphi\rangle\geq (\sigma_{p+1}(M)-g(S(\xi),\xi))|\varphi|^2,$$
for any basic $p$-form $\varphi$ on $M.$ 
\end{lemma}

Moreover, we need to get an upper bound for the norm of the interior product of the $2$-form $\Omega$ with any basic $p$-form. Indeed,

\begin{lemma}\label{lem:4}
Let $(N^{n+2},g)$ be a Riemannian manifold with boundary $M$. Assume that $M$ carries a Riemannian flow given by a unit vector field $\xi$. For any basic $p$-form $\varphi$ with $p\geq 2,$ we have
$$|\Omega\lrcorner\varphi|\leq [\frac{n}{2}]^{\frac{1}{2}}|\Omega| |\varphi|.$$
If $n$ is even, the equality is realized if and only if $\Omega=0$.
\end{lemma}
{\bf Proof.} Let $\lambda_j$ be the eigenvalues of the $2$-form $\Omega$. We can always find an orthonormal frame $\{e_i\}$ of $\Gamma(Q)$ such that $\Omega=\sum_{j=1}^{[\frac{n}{2}]}\lambda_j e_{2j-1}\wedge e_{2j}$. Therefore, we compute
\begin{eqnarray*}
|\Omega\lrcorner\varphi|&=&|\sum_{j=1}^{[\frac{n}{2}]}\lambda_j (e_{2j-1}\wedge e_{2j})\lrcorner \varphi|\\
&\leq & \sum_{j=1}^{[\frac{n}{2}]}|\lambda_j| |e_{2j-1}\lrcorner(e_{2j}\lrcorner\varphi)|\\
&\leq & \sum_{j=1}^{[\frac{n}{2}]}|\lambda_j||\varphi|\leq [\frac{n}{2}]^{\frac{1}{2}}|\Omega| |\varphi|.
\end{eqnarray*}
Here we used the fact that $|v\lrcorner \varphi|\leq |v||\varphi|$ and the Cauchy-Schwarz inequality (in the last estimate). Assume now that the equality is realized, then either all the $\lambda_j's$ are of the same absolute value and there exists a $j$ such that $\lambda_j=0$ (in this case, all the $\lambda_j$'s are 0) or for all $j$, $e_{2j}\wedge \varphi=0$ and $e_{2j-1}\wedge \varphi=0.$ But for $n$ even, the last statement just means that $X\wedge \varphi=0$ for all $X\in \Gamma(Q)$ and thus $\varphi=0.$ This leads to a contradiction; hence $\lambda_j=0$ for all $j$ which yields $\Omega=0$.
\hfill$\square$\\


Now, we have all the materials to prove Theorem \ref{thm:main}:\\
 
{\noindent \bf Proof of Theorem \ref{thm:main}.} For any basic closed $p$-eigenform  $\varphi$ corresponding to the eigenvalue $\lambda'_{1,p}$ of the basic Laplacian, we associate its extension $\hat\varphi$ that is closed and co-closed on $N$  from Lemma \ref{lem:1}. 
Applying the Reilly formula to the $p$-form $\hat\varphi$ gives, under the curvature assumption and with the use of Lemma \ref{lem:2} for the eigenform $\varphi$, the following
$$
0\geq 2\int_M\langle \nu\lrcorner\hat\varphi, \delta^M\varphi\rangle+ \sigma_{p+1}(M)\int_M|\varphi|^2-\int_Mg(S(\xi),\xi)|\varphi|^2\\+\sigma_{n+2-p}(M)\int_M|\nu\lrcorner\hat\varphi|^2.$$
Using the pointwise inequality $|\nu\lrcorner\hat\varphi+\frac{1}{\sigma_{n+2-p}(M)}\delta^M\varphi|^2\geq 0,$ the above one can be reduced to the following
\begin{equation}\label{eq:25}
\int_M|\delta^M\varphi|^2\geq \sigma_{p+1}(M) \sigma_{n+2-p}(M)\int_M|\varphi|^2- \sigma_{n+2-p}(M)\mathop{{\rm sup}}\limits_M(g(S(\xi),\xi))\int_M|\varphi|^2.
\end{equation}
Now from the relation $\delta_b=\delta_M-2\Omega\lrcorner(\xi\wedge)$ on basic forms \cite[Prop.2.4]{RP} and the estimate in Lemma \ref{lem:4}, we get
\begin{eqnarray*}
|\delta^M\varphi|^2&=& |\delta_b\varphi|^2+4|\Omega\lrcorner\varphi|^2+4\langle\delta_b\varphi,\Omega\lrcorner(\xi\wedge\varphi)\rangle\\
&\leq & |\delta_b\varphi|^2+4c^2 [\frac{n}{2}]|\varphi|^2.
\end{eqnarray*}
The third term is zero, since $\delta_b\varphi$ is basic. Therefore after integrating over the manifold $M,$ we finish the proof of the theorem.
\hfill$\square$ 

\begin{remark}
The assumptions on the curvature could be weakened. The positivity of the curvature operator could be replaced by the positivity of $W_N^{[p]}$ and $W_N^{[p+1]}$.
\end{remark}

\section{The equality case}
This section is devoted to the proof of Theorem $\ref{thm:mainestimate}$. In other words, we are going to study 
the limiting case of Inequality \eqref{ineq:main}. We will show that, under some conditions, the second fundamental form vanishes along $\xi$ and is equal to $\eta\, {\rm Id}$ in the direction of $Q$ for some constant $\eta,$ i.e. the boundary is $\eta$-umbilical.
 We will also prove that the O'Neill tensor defining the flow vanishes; this is equivalent to the integrability of the normal bundle. 
 Consequently, the extension of the vector field $\xi$ given by Problem \eqref{eq:23} is parallel on the whole manifold. This allows the classification of all manifolds on which Inequality \eqref{ineq:main} is optimal.

It is clear to see that when the equality is realized, the estimate in Lemma \ref{lem:4} is optimal which means that $h=0$. On the other hand, the eigenform $\hat\varphi$ is parallel on $N$ and $\sigma_{p+1}, \sigma_{n+2-p}$ and $g(S(\xi),\xi)$ are constant on $M$. Moreover, we have the relation
\begin{equation}\label{eq:26}
\delta^M\varphi=-\sigma_{n+2-p}\nu\lrcorner \hat\varphi.
\end{equation}
In particular, using the relations in \cite[Lemma 18]{RS}, we get for all $X\in \Gamma(TM)$,
\begin{equation} \label{eqcase}
\left\{
\begin{array}{lll}
\nabla^M_X \varphi=S(X)\wedge (\nu\lrcorner \hat\varphi) \\\\ \nabla^M_X(\nu\lrcorner \hat\varphi)=-S(X) \lrcorner\varphi,\\\\
\delta^M \varphi =S^{[p-1]}(\nu\lrcorner\hat\varphi)-\sigma_{n+1} \nu\lrcorner \hat\varphi \\\\
d^M (\nu\lrcorner \hat\varphi)=-S^{[p]}(\varphi).
\end{array}
\right.
\end{equation}
Hence using the last two equations in \eqref{eqcase} along with Equation \eqref{eq:26}, we deduce that
\begin{equation}\label{eq:27}
\left\{
\begin{array}{lll}
S^{[p-1]}(\nu\lrcorner\hat\varphi)=(\sigma_{n+1}-\sigma_{n+2-p})\nu\lrcorner\hat\varphi\\\\
S^{[p]}(\varphi)=(\sigma_{p+1}-g(S(\xi),\xi))\varphi.
\end{array}
\right.
\end{equation}
In order to prove that the manifold $N$ is isometric to the quotient  $\lquot{\mathbb{R} \times B'}{\Gamma}$, we need first to prove a series of lemmas:

\begin{lemma}
If the equality is realized in $(\ref{ineq:main})$, then $S(\xi)=0$.
\end{lemma}
{\bf Proof.} 
 Using Equation \eqref{eq:26}, we deduce that the form $\nu\lrcorner \hat\varphi$ is basic (recall here that $\delta^M\varphi=\delta_b\varphi$). Hence by applying the first equation in \eqref{eq:27} to the vector fields $\xi$ and $X_1,\cdots,X_{p-2} \in \Gamma(Q)$, we find that $S(\xi)\lrcorner(\nu\lrcorner\hat\varphi)=0$. On the other hand, since the O'Neill tensor vanishes, then $\nabla^M_\xi\varphi=\nabla_\xi\varphi$ which is equal to zero, because the form $\varphi$ is basic. Here, we recall that $\nabla$ is the extension of the transversal Levi-Civita connection $\nabla$ to basic forms. Finally, by taking $X=\xi$ in the first equation of \eqref{eqcase} we find that $S(\xi)\wedge (\nu\lrcorner\hat\varphi)=0$. Mainly, that means $S(\xi)=0$. We mention here that $\nu\lrcorner\hat\varphi$ cannot vanish, since this would imply that $\nabla^M_X\varphi=0$ for all $X\in \Gamma(TM)$ which would give that $\lambda'_{1,p}=0.$
\hfill$\square$ \\ 

In the sequel, we aim to prove that the principal curvatures of $S$ are constant and are all equal to a number $\eta$, along transversal principal directions. The proof of this statement is a technical computation and will be splitted into several lemmas (see Lemmas \ref{lem:geo1}, \ref{lem:geo2} and \ref{lem:geo3}). For this, we will consider in all our calculations an orthonormal frame $\{f_i\}_{i=1,\cdots,n+1}$ of $\Gamma(TM)$.

\begin{lemma} \label{lem:geo1}
If  the equality is realized in $(\ref{ineq:main}),$ the identity
\begin{eqnarray} \label{eq:chara}
\sum_{i=1}^{n+1}\langle (\nabla^M_{f_i} S)^{[p]}\varphi,f_i\wedge (\nu\lrcorner\hat\varphi)\rangle&=&((\sigma_{p+1}-\sigma_{n+1}+\sigma_{n+2-p})\sigma_{n+2-p}-|S|^2)|\nu\lrcorner\hat \varphi|^2\nonumber\\
&&+\sum_{i=1}^{n+1}\langle f_i\lrcorner(\nu\lrcorner\hat\varphi),S^2(f_i)\lrcorner (\nu\lrcorner\hat\varphi)\rangle
\end{eqnarray}
holds.
\end{lemma}
{\bf Proof.}
By differentiating the second equation in \eqref{eq:27} in the direction of any vector field $X\in \Gamma(TM)$, we get after using \eqref{eq: appe1}
$$S^{[p]}(SX\wedge (\nu\lrcorner \hat\varphi))+(\nabla^M_X S)^{[p]}\varphi=\sigma_{p+1}SX\wedge (\nu\lrcorner \hat\varphi).$$
Here we also used the first equation in \eqref{eqcase}.
Setting $X=f_i$ and taking the scalar product of the last equality with $f_i\wedge (\nu\lrcorner\hat\varphi)$, we obtain after tracing and using \eqref{eq: appe2} that,
\begin{eqnarray*}
\sum_{i=1}^{n+1}\langle (\nabla^M_{f_i} S)^{[p]}\varphi,f_i\wedge (\nu\lrcorner\hat\varphi)\rangle&=&\sigma_{p+1}\sum_{i=1}^{n+1}\langle S(f_i)\wedge (\nu\lrcorner \hat\varphi),f_i\wedge (\nu\lrcorner\hat\varphi)\rangle\nonumber\\
&&-\sum_{i=1}^{n+1}\langle S^2 (f_i)\wedge(\nu\lrcorner\hat\varphi),f_i\wedge (\nu\lrcorner\hat\varphi)\rangle\nonumber\\
&&-\sum_{i=1}^{n+1}\langle S (f_i) \wedge S^{[p-1]}(\nu\lrcorner \hat\varphi),f_i\wedge(\nu\lrcorner \hat\varphi)\rangle. 
\end{eqnarray*} 
Then with the help of \eqref{eq:27}, the last equality reduces to 
\begin{eqnarray}\label{eq:trace}
\sum_{i=1}^{n+1}\langle (\nabla^M_{f_i} S)^{[p]}\varphi,f_i\wedge (\nu\lrcorner\hat\varphi)\rangle
&=&(\sigma_{p+1}-\sigma_{n+1} +\sigma_{n+2-p})\sum_{i=1}^{n+1}\langle S(f_i)\wedge (\nu\lrcorner \hat\varphi),f_i\wedge (\nu\lrcorner\hat\varphi)\rangle \nonumber\\ &&-\sum_{i=1}^{n+1}\langle S^2(f_i)\wedge(\nu\lrcorner\hat\varphi),f_i\wedge (\nu\lrcorner\hat\varphi)\rangle.
\end{eqnarray}
In order to finish the proof, it is sufficient to calculate the two sums in the r.h.s. of \eqref{eq:trace}. In fact, the first sum is equal to
\begin{eqnarray*}
\sum_{i=1}^{n+1}\langle S(f_i)\wedge (\nu\lrcorner \hat\varphi),f_i\wedge (\nu\lrcorner\hat\varphi)\rangle&=
&\sigma_{n+1}|\nu\lrcorner \hat\varphi|^2-\sum_{i=1}^{n+1}\langle f_i\lrcorner (\nu\lrcorner\hat\varphi), S(f_i)\lrcorner (\nu\lrcorner\hat\varphi)\rangle\\
&\bui{=}{\eqref{eq:appe4},\eqref{eq:27}}&\sigma_{n+2-p}|\nu\lrcorner \hat\varphi|^2,
\end{eqnarray*} 
while the second one is
\begin{equation*}
\sum_{i=1}^{n+1}\langle S^2(f_i)\wedge(\nu\lrcorner\hat\varphi),f_i\wedge (\nu\lrcorner\hat\varphi)\rangle=
|S|^2|\nu\lrcorner\hat \varphi|^2-\sum_{i=1}^{n+1}\langle S^2 (f_i) \lrcorner(\nu\lrcorner\hat\varphi),f_i\lrcorner (\nu\lrcorner\hat\varphi)\rangle.
\end{equation*}
The substitution into \eqref{eq:trace} gives the desired result.
\hfill$\square$\\ 

\noindent In the following lemma, we will compute the l.h.s. of Equation \eqref{eq:chara} in terms of the curvature operator $W^{[p]}$. Indeed,

\begin{lemma}\label{lem:geo2}
If the equality is realized  in $(\ref{ineq:main})$, the relation
\begin{eqnarray} \label{eq:left}
\sum_{i=1}^{n+1}\langle (\nabla^M_{f_i} S)^{[p]}\varphi,f_i\wedge (\nu\lrcorner\hat\varphi)\rangle&=&-\langle W^{[p]}(\nu\wedge(\nu\lrcorner\hat\varphi)),\varphi\rangle+\sum_{i=1}^{n+1}\langle f_i\lrcorner\varphi,S^2(f_i)\lrcorner\varphi\rangle\nonumber\\
&&+(\sigma_{n+1}-\sigma_{n+2-p}-\sigma_{p+1})\sigma_{p+1}|\varphi|^2,
\end{eqnarray}
holds.
\end{lemma}
{\bf Proof.} Using the symmetry property of the tensor $\nabla^M S$ and Equation \eqref{eq: appe2}, the l.h.s. of Equation \eqref{eq:chara} is equal to
\begin{eqnarray*}
\sum_{i=1}^{n+1}\langle \varphi,(\nabla^M_{f_i} S)(f_i)\wedge(\nu\lrcorner\hat\varphi)\rangle+\langle \varphi,f_i\wedge(\nabla^M_{f_i}S)^{[p-1]}(\nu\lrcorner\hat\varphi)\rangle\\
\bui{=}{\eqref{eq:appe3}}(-1)^{\frac{p(p-1)}{2}}\sum_{i=1}^{n+1}\langle(\nu\lrcorner\hat\varphi)\lrcorner\varphi,(\nabla^M_{f_i} S)(f_i)\rangle
+\sum_{i=1}^{n+1}\langle \varphi,f_i\wedge(\nabla^M_{f_i}S)^{[p-1]}(\nu\lrcorner\hat\varphi)\rangle.
\end{eqnarray*} 
Therefore, from Equation \eqref{eq:codazzi} and again from the symmetry of $\nabla^M S,$ the above expression reduces to 
\begin{eqnarray} \label{eqcase1}
(-1)^{\frac{p(p-1)}{2}}\sum_{i=1}^{n+1}\langle(\nabla^M_{(\nu\lrcorner\hat\varphi)\lrcorner\varphi}S)(f_i),f_i\rangle
+(-1)^{\frac{p(p-1)}{2}}R^N((\nu\lrcorner\hat\varphi)\lrcorner\varphi,f_i,\nu,f_i)
+\sum_{i=1}^{n+1}\langle \varphi,f_i\wedge(\nabla^M_{f_i}S)^{[p-1]}(\nu\lrcorner\hat\varphi)\rangle\nonumber\\
=(-1)^{\frac{p(p-1)}{2}}((\nu\lrcorner\hat\varphi)\lrcorner\varphi)(\sigma_{n+1})
-(-1)^{\frac{p(p-1)}{2}}\sum_{i=1}^{n+1}R^N(\nu,f_i,f_i,(\nu\lrcorner\hat\varphi)\lrcorner\varphi)
+\sum_{i=1}^{n+1}\langle \varphi,f_i\wedge(\nabla^M_{f_i}S)^{[p-1]}(\nu\lrcorner\hat\varphi)\rangle.\nonumber\\
\end{eqnarray}
Let us explicit now the curvature term $\langle W^{[p]}(\nu\wedge(\nu\lrcorner\hat\varphi)),\varphi\rangle.$ Recall that $W^{[p]}=\sum_{i,j} e_j^*\wedge e_i\lrcorner R^N(e_i,e_j)$ where $\{e_i\}_{i=1,\cdots,n+2}$ is any orthonormal frame of $TN$. At a point $x\in M$, we take the orthonormal frame on $T_x N$ as $\{f_i,\nu\}_{i=1,\cdots,n+1}$ and get

\begin{eqnarray*}
\langle W^{[p]}(\nu\wedge(\nu\lrcorner\hat\varphi)),\varphi\rangle&=&\sum_{i=1}^{n+1}\langle \nu\lrcorner R^N(\nu,f_i)(\nu\wedge(\nu\lrcorner\hat\varphi)), f_i\lrcorner\varphi\rangle
+\sum_{i,j=1}^{n+1}\langle f_i\lrcorner R^N(f_i,f_j)(\nu\wedge(\nu\lrcorner\hat\varphi)), f_j\lrcorner\varphi\rangle\nonumber\\
&=&\sum_{i=1}^{n+1}\langle R^N(\nu,f_i)\nu\lrcorner\hat\varphi, f_i\lrcorner\varphi\rangle
+\sum_{i,j=1}^{n+1}\langle R^N(\nu,f_i)f_i,f_j\rangle\langle \nu\lrcorner\hat\varphi,f_j\lrcorner\varphi\rangle
\nonumber\\&&-\sum_{i,j=1}^{n+1}\langle R^N(f_i,f_j)\nu\wedge f_i\lrcorner(\nu\lrcorner\hat\varphi),f_j\lrcorner\varphi\rangle.
\end{eqnarray*} Then using \eqref{eq:appe3}, we deduce that
\begin{eqnarray}\label{cur}
\langle W^{[p]}(\nu\wedge(\nu\lrcorner\hat\varphi)),\varphi\rangle& =&\sum_{i=1}^{n+1}\langle R^N(\nu,f_i)\nu\lrcorner\hat\varphi, f_i\lrcorner\varphi\rangle
+(-1)^{\frac{p(p-1)}{2}}\sum_{i=1}^{n+1} R^N(\nu,f_i,f_i,(\nu\lrcorner\hat\varphi)\lrcorner\varphi)
\nonumber\\
&&-\sum_{i,j=1}^{n+1}\langle R^N(f_i,f_j)\nu\wedge f_i\lrcorner(\nu\lrcorner\hat\varphi),f_j\lrcorner\varphi\rangle\nonumber.\\
\end{eqnarray}
For simplicity, we will denote the first, the second and the third terms in \eqref{cur} respectively by $(*), (**)$ and $(***).$ On one hand, we remark that the term $(***)$ is equal to the following:
\begin{eqnarray} \label{eq:337}
\sum_{\substack{i_1<\cdots<i_{p-1}\\k=1\cdots,p-1\\i,j=1,\cdots,n+1}}&(-1)^{k+1}(\nu\lrcorner\hat\varphi)_{i_1,\cdots,i_{p-1}}\delta_{ii_k}\langle R^N(f_i,f_j)\nu\wedge f_{i_1}\wedge \cdots {\hat f_{i_k}}\wedge \cdots\wedge f_{i_{p-1}},f_j\lrcorner\varphi\rangle \nonumber\\
&=\displaystyle \sum_{\substack{i_1<\cdots<i_{p-1}\\k=1\cdots,p-1\\j=1,\cdots,n+1}}(\nu\lrcorner\hat\varphi)_{i_1,\cdots,i_{p-1}}\langle  f_{i_1}\wedge \cdots \wedge R^N(f_{i_k},f_j)\nu\wedge \cdots\wedge f_{i_{p-1}},f_j\lrcorner\varphi\rangle \nonumber\\
&\bui{=}{\eqref{eq:codazzi}}\displaystyle \sum_{\substack{i_1<\cdots<i_{p-1}\\k= 1\cdots,p-1\\j=1,\cdots,n+1}}(\nu\lrcorner\hat\varphi)_{i_1,\cdots,i_{p-1}}\langle  f_{i_1}\wedge \cdots \wedge (\nabla^M_{f_j}S)(f_{i_k})\wedge \cdots\wedge f_{i_{p-1}},f_j\lrcorner\varphi\rangle \nonumber\\
&=\displaystyle \sum_{j=1}^{n+1}\langle \varphi,f_j\wedge(\nabla^M_{f_j}S)^{[p-1]}(\nu\lrcorner\hat\varphi)\rangle\nonumber.\\
\end{eqnarray}
Hence Equation \eqref{eqcase1} becomes
\begin{eqnarray}\label{eq:336}
\sum_{i=1}^{n+1}\langle (\nabla^M_{f_i} S)^{[p]}\varphi,f_i\wedge (\nu\lrcorner\hat\varphi)\rangle&=&(-1)^{\frac{p(p-1)}{2}}((\nu\lrcorner\hat\varphi)\lrcorner\varphi)(\sigma_{n+1})-\langle W^{[p]}(\nu\wedge(\nu\lrcorner\hat\varphi)),\varphi\rangle\nonumber
\\&&+(*).
\end{eqnarray}
On the other hand, one can easily check that $(*)+(***)=0$. Indeed, using the second equation in \eqref{eq:337}, the sum of $(*)$ and $(***)$ is equal to the following
\begin{eqnarray*}
\sum_{\substack{i_1<\cdots<i_{p-1}\\k=1\cdots,p-1\\i=1,\cdots,n+1}}(\nu\lrcorner\hat\varphi)_{i_1,\cdots,i_{p-1}}\{\langle f_{i_1}\wedge\cdots\wedge R^N(\nu,f_i)f_{i_k}\wedge\cdots \wedge f_{i_{p-1}},f_i\lrcorner \varphi\rangle\\
+\langle  f_{i_1}\wedge \cdots \wedge R^N(f_{i_k},f_i)\nu\wedge \cdots\wedge f_{i_{p-1}},f_i\lrcorner\varphi\rangle\}\\
=\sum_{\substack{i_1<\cdots<i_{p-1}\\k=1,\cdots,p-1\\i,l=1,\cdots,n+1}}(\nu\lrcorner\hat\varphi)_{i_1,\cdots,i_{p-1}}\{R^N(f_{i_k},f_l,\nu,f_i)\langle f_{i_1}\wedge\cdots\wedge f_l\wedge\cdots \wedge f_{i_{p-1}},f_i\lrcorner \varphi\rangle\\
+R^N(f_{i_k},f_i,\nu,f_l)\langle  f_{i_1}\wedge \cdots \wedge f_l\wedge \cdots\wedge f_{i_{p-1}},f_i\lrcorner\varphi\rangle\},
\end{eqnarray*}
which is zero when one interchanges the role of the indices $i$ and $l$ in the first summation. Therefore, and since $(*)=-(***)=\displaystyle -\sum_{i=1}^{n+1}\langle f_i\lrcorner\varphi,(\nabla^M_{f_i}S)^{[p-1]}(\nu\lrcorner\hat\varphi)\rangle $, then 
\begin{eqnarray*}
(*)
&\bui{=}{\eqref{eq: appe1}}&-\sum_{i=1}^{n+1}\langle f_i\lrcorner\varphi,\nabla^M_{f_i}(S^{[p-1]}(\nu\lrcorner\hat\varphi))\rangle+\sum_{i=1}^{n+1}\langle f_i\lrcorner\varphi,S^{[p-1]}(\nabla^M_{f_i}(\nu\lrcorner\hat\varphi))\rangle\nonumber\\
&\bui{=}{\eqref{eqcase}, \eqref{eq:27}}&-\sum_{i=1}^{n+1}f_i(\sigma_{n+1})\langle f_i\lrcorner\varphi,\nu\lrcorner\hat\varphi\rangle+(\sigma_{n+1}-\sigma_{n+2-p})\sum_{i=1}^{n+1}\langle f_i\lrcorner\varphi,S(f_i)\lrcorner\varphi\rangle\nonumber\\
&&-\sum_{i=1}^{n+1}\langle f_i\lrcorner\varphi,S^{[p-1]}(S(f_i)\lrcorner\varphi)\rangle, \end{eqnarray*}
which is, using Equations \eqref{eq:appe3},\eqref{eq:appe4},\eqref{eq:appe5} and \eqref{eq:27}
\begin{eqnarray}\label{eq:338}
-(-1)^{\frac{p(p-1)}{2}}((\nu\lrcorner\hat\varphi)\lrcorner\varphi)(\sigma_{n+1})+(\sigma_{n+1}-\sigma_{n+2-p})\sigma_{p+1}|\varphi|^2
-\sigma_{p+1}^2|\varphi|^2+\sum_{i=1}^{n+1}\langle f_i\lrcorner\varphi,S^2(f_i)\lrcorner\varphi\rangle.
\end{eqnarray}
Substituting Equation \eqref{eq:338} into Equation \eqref{eq:336}, we finally get the result.
\hfill$\square$ \\

\noindent In the next lemma, we shall compare the sign of the l.h.s. of Equation \eqref{eq:chara} which is given by \eqref{eq:left} to the r.h.s. and shall find that they are of opposite signs, when a curvature assumption is required. In particular, this will mean that all principal curvatures along transversal directions are equal. The statement is:

\begin{lemma}\label{lem:geo3}
If the equality is realized and if moreover $\sigma_1(M)\geq 0$, then $S(X)=\eta X$ for all $X\in \Gamma(Q).$
\end{lemma}
{\bf Proof.} We will show that the l.h.s. of Equation \eqref{eq:chara} is non-negative while the r.h.s. is non-positive. We first begin to check the l.h.s. The eigenform $\hat \varphi$ being parallel, the term $\langle W^{[p]}\hat\varphi,\hat\varphi\rangle$ vanishes. That mainly means, by writing $\hat\varphi=\varphi+\nu\wedge(\nu\lrcorner\hat\varphi)$ at any point of the boundary and using the fact that $W^{[p]}$ is non-negative, the term $\langle W^{[p]}(\nu\wedge(\nu\lrcorner\hat\varphi)),\varphi\rangle$ is non-positive. On the other hand,  the tensor $S$ has $0$ as an eigenvalue (recall that $S(\xi)=0$) and that $\sigma_1(M)\geq 0$, then all $\eta_i$'s are greater than $0$ for $i=2,\cdots,n+1.$ Hence, we get the estimate
$$
\sum_{i=1}^{n+1}\langle f_i\lrcorner\varphi,S^2(f_i)\lrcorner\varphi\rangle\geq \eta_2\sum_{i=1}^{n+1}\langle f_i\lrcorner\varphi,S(f_i)\lrcorner\varphi\rangle\bui{=}{\eqref{eq:appe4},\eqref{eq:27}}\eta_2\sigma_{p+1}|\varphi|^2.
$$
Therefore, Equation \eqref{eq:left} allows to bound from below the l.h.s. by
\begin{eqnarray*}
&(\sigma_{n+1}-\sigma_{n+2-p}-\sigma_{p+1}+\eta_2)\sigma_{p+1}|\varphi|^2\\
=&((\eta_{n+3-p}-\eta_3) +\cdots +(\eta_{n+1}-\eta_{p+1}))\sigma_{p+1}|\varphi|^2
\geq 0,
\end{eqnarray*}
since the sequence $\eta_i$ is increasing. We can easily see that the last expression vanishes when all the $\eta_i's$ are equal. Concerning the r.h.s. of Equation \eqref{eq:chara}, recall that it is given by
\begin{equation}\label{right}
((\sigma_{p+1}-\sigma_{n+1}+\sigma_{n+2-p})\sigma_{n+2-p}-|S|^2)|\nu\lrcorner\hat \varphi|^2
+\sum_{i=1}^{n+1}\langle f_i\lrcorner(\nu\lrcorner\hat\varphi),S^2(f_i)\lrcorner (\nu\lrcorner\hat\varphi)\rangle.
\end{equation}
In the sequel, we will take the vectors $\{f_i\}_{i=1,\cdots,n+1}$ as the principal directions associated with the principal curvatures $\eta_i$ of the tensor $S$. We first estimate
\begin{eqnarray*}
\sum_{i=1}^{n+1}\langle f_i\lrcorner(\nu\lrcorner\hat\varphi),S^2(f_i)\lrcorner (\nu\lrcorner\hat\varphi)\rangle &=&\sum_{i=1}^{n+1}\eta_i\langle f_i\lrcorner(\nu\lrcorner\hat\varphi),S(f_i)\lrcorner (\nu\lrcorner\hat\varphi)\rangle\\
&\leq &\eta_{n+1}\sum_{i=1}^{n+1}\langle f_i\lrcorner(\nu\lrcorner\hat\varphi),S(f_i)\lrcorner (\nu\lrcorner\hat\varphi)\rangle\\
&\bui{=}{\eqref{eq:appe4},\eqref{eq:27}}&\eta_{n+1}(\sigma_{n+1}-\sigma_{n+2-p})|\nu\lrcorner\hat\varphi|^2.
\end{eqnarray*}
Hence, \eqref{right} can be bounded from above by
$$\eqref{right} \le A |\nu\lrcorner\hat\varphi|^2,$$
where $A$ is given by
$$A = (\sigma_p-\sigma_{n+1}+\sigma_{n+2-p})\sigma_{n+2-p}+\eta_{p+1}\sigma_{n+2-p}-\eta_2^2-\cdots-\eta_n^2
+\eta_{n+1}(\eta_{n+3-p}+\cdots+\eta_n).$$
In the following, we will prove that $A$ is non-positive, which implies that \eqref{right} is non-positive.
\begin{eqnarray*}
A&= & (\sigma_p-\sigma_{n+1}+\sigma_{n+2-p})\sigma_{p+1}+(\sigma_p-\sigma_{n+1}+\sigma_{n+2-p})(\eta_{p+2}+\cdots+\eta_{n+2-p})
\\
&&+\eta_{p+1}\sigma_{p}+\eta_{p+1}(\eta_{p+1}+\cdots+\eta_{n+2-p})-\eta_2^2-\cdots-\eta_n^2
+\eta_{n+1}(\eta_{n+3-p}+\cdots+\eta_n)\\
&= &\sigma_p \sigma_{p+1}-(\sigma_{n+1}-\sigma_{n+2-p})\sigma_{p+1}+(\sigma_p-\sigma_{n+1}+\sigma_{n+2-p})(\eta_{p+2}+\cdots+\eta_{n+2-p})
\\
&&+\eta_{p+1}\sigma_{p} +B-\eta_2^2-\cdots-\eta_p^2-\eta_{n+3-p}^2-\cdots-\eta_n^2+\eta_{n+1}(\eta_{n+3-p}+\cdots+\eta_n)
\end{eqnarray*}
where $B$ is given  by
$$B=\eta_{p+2}(\eta_{p+1}-\eta_{p+2})+\eta_{p+3}(\eta_{p+1}-\eta_{p+3})+\cdots+\eta_{n+2-p}(\eta_{p+1}-\eta_{n+2-p}).$$
\noindent
Clearly $B$ is non-positive. One the other hand, since
$$
S^{[p]}(\xi\wedge(\nu\lrcorner\hat\varphi))\bui{=}{\eqref{eq: appe2}}\xi\wedge S^{[p-1]}(\nu\lrcorner\hat\varphi)\bui{=}{\eqref{eq:27}}(\sigma_{n+1}-\sigma_{n+2-p})(\xi\wedge(\nu\lrcorner\hat\varphi)),$$ then $\sigma_p+\sigma_{n+2-p}\leq \sigma_{n+1}.$
Here we use the fact that $S(\xi)=0.$ Therefore 
\begin{eqnarray*}
A &\leq  &\sigma_p^2+2\eta_{p+1}\sigma_{p}-(\sigma_{n+1}-\sigma_{n+2-p})\sigma_{p+1}-\eta_2^2-\cdots-\eta_p^2\\&&-\eta_{n+3-p}^2-\cdots-\eta_n^2
+\eta_{n+1}(\eta_{n+3-p}+\cdots+\eta_n)\\
&=&2\sum_{2\leq i<j\leq p}\eta_i\eta_j+2\eta_{p+1}\sigma_{p}-\eta_{n+3-p}(\eta_{n+3-p}+\sigma_{p+1})-\cdots-\eta_n(\eta_n+\sigma_{p+1})\\
&&-\eta_{n+1}(\sigma_{p+1}-\eta_{n+3-p}-\cdots-\eta_n).
\end{eqnarray*}
Using the fact that for $i=n+3-p, \cdots, n+1$ each $\eta_i\geq \eta_{p+2}$ and that $\eta_{n+3-p}+\cdots+\eta_n\leq \sigma_{p+1},$ since $\sigma_{p+1}-(\eta_{n+3-p}+\cdots+\eta_n)$ is an eigenvalue of $S^{[2]}$ (just apply $S^{[2]}$ to the eigenform $(f_{i_{n+3-p}}\wedge\cdots\wedge f_{i_{n}})\lrcorner\varphi$ by using the formula \eqref{eq:appe6}), we deduce that
$$A \leq 2\sum_{2\leq i<j\leq p}\eta_i\eta_j+2\eta_{p+1}\sigma_{p}-\eta_{p+2}(p-1)\sigma_{p+1}\leq 0.
$$
This last inequality is true because the number of positive terms is equal to the number of negative terms which is $p(p-1).$
\hfill$\square$\\

Now we are able to prove Theorem  \ref{thm:mainestimate}.\\

{\noindent\bf Proof of Theorem \ref{thm:mainestimate}.} We proceed as in \cite{EHGN1} (see also \cite{EHGN} for more details). We first show that the vector field $\xi$
defining the flow can be extended to a unique parallel vector field
$\hat\xi$ on $N$ which is orthogonal to $\nu$. The proof mainly relies on the use of the Reilly formula on the solution $\hat\xi$ of the boundary problem \eqref{eq:23}. Second, we consider a connected integral submanifold $N_1$ of the bundle $(\R\hat{\xi})^\perp,$ where the orthogonal is taken in $N$. The manifold $N_1$ is complete with totally umbilical boundary and the Ricci tensor of $\partial N_1$ is bounded from below by some constant. That means, the manifold $\partial N_1$ is compact as a consequence of Myers's theorem. This allows to deduce that $N_1$ is compact from the main theorem in \cite[Thm. 1.1]{Li}.

\noindent On the other hand, we have from Equations \eqref{eqcase} that $\varphi=-\frac{1}{p\eta}d^M(\nu\lrcorner \hat\varphi)$
which means that it is $d^M$-exact and thus $d^{\partial N_1}$-exact, since $\partial N_1$ is totally geodesic in $M$ and
both $\varphi$ and  $\nu\lrcorner \hat\varphi$ are basic. Moreover, the basic form $\varphi$ is an eigenform of the Laplacian on $\partial N_1$, that is $\Delta^{\partial N_1}\varphi=\Delta_b\varphi=\lambda'_{1,p}\varphi.$ Therefore, if we denote by $\lambda_{1,p}^{\partial N_1}$ the first eigenvalue
of $\Delta^{\partial N_1}$ restricted to exact $p$-forms on $N_1$ and by $\tilde\sigma_p$ the $p$-curvatures of $\partial N_1$
into the compact manifold $N_1,$ we get from the main estimate in \cite[Thm. 5]{RS} that
$$p(n+1-p)\eta^2=\tilde\sigma_p \tilde\sigma_{n+1-p}\leq \lambda_{1,p}^{\partial N_1}\leq \lambda'_{1,p}=\sigma_{p+1}\sigma_{n+2-p}=p(n+1-p)\eta^2.$$
Hence the equality is attained in the estimate of Raulot and Savo and therefore $N_1$ is isometric to the Euclidean closed ball $B'$. Finally, by the de Rham theorem, the manifold $\widetilde{N}$ is
isometric to $\mathbb{R}\times B'$ and $N$ is the
quotient of the Riemannian product $\R\times B'$ by its
fundamental group. Since $\pi_1(N)$ embeds into $\pi_1(M)$ 
, $N$ is isometric to $\lquot{\mathbb{R} \times
B'}{\Gamma}$.
\hfill$\square$

\section{Rigidity results on manifolds with foliated boundary}\label{sec:rigidity}
Our objective, in this section, is to derive rigidity results on manifolds with foliated boundary. These results generalize the ones in \cite[Sect. 5]{EHGN1}. For this end, we recall that a basic special Killing $p$-form $\omega$ is a basic co-closed (with respect to the basic codifferential $\delta_b$) form satisfying for all $X \in \Gamma(Q)$ the relations
$$\nabla_X\omega = \frac{1}{p+1} X\lrcorner d_b\omega \ \ \ \ \ \  \text{and} \ \ \ \ \nabla_X d_b\omega = -c(p+1) X \wedge \omega,$$
where $\nabla$ is the transversal Levi-Civita connection, defined in Section 2, extended to basic forms and $c$ is a non-negative constant. In general, one can prove that a basic special Killing $p$-form is a co-closed eigenform of the basic Laplacian corresponding to the eigenvalue $c(p+1)(n-p)$ where $n$ is the rank of $Q$.

\noindent In the following, we will consider a compact manifold $N$ whose boundary carries a basic special Killing $p$-form. We will see how we could characterize the boundary as the product $\mathbb{S}^1\times \mathbb{S}^n$ and this is due to the equality case of our main estimate. 
\noindent We first prove the following result:

\begin{cor}\label{equalitycara}
Let $N$ be an $(n+2)$-dimensional compact manifold with non-negative
curvature operator. Assume that the boundary $M$ carries a minimal
Riemannian flow such that $(n+1-p)\sup_M g(S(\xi),\xi) + 4c^2 [\frac n2]\leq 0$ and also admits a
basic special Killing $(n-p)$-form for some $2 \leq p \leq \frac n2$.  If the inequality
$\sigma_{p+1}(M)\geq p$ holds, the manifold $N$ is isometric to
$\lquot{\mathbb{R} \times B'}{\Gamma}$.
\end{cor}

{\bf Proof.} Let $\varphi$ be a basic special Killing $(n-p)$-form on $M$. Then $*_b\varphi$ is a basic closed $p$-eigenform for the basic Laplacian, that is $\Delta_b(*_b\varphi)=p (n+1-p)(*_b\varphi).$ Here, we used the minimality of the flow to say that the basic Hodge operator ``$*_b$'' commutes with the basic Laplacian. Hence $\lambda'_{1,p}\leq p(n+1-p).$ To get the upper bound, we will use the estimate in Theorem \ref{thm:main}. First, we have $\sigma_{n+1-p} \geq \sigma_{p+1} >0$ as $2 \leq p \leq \frac n2.$
On the other hand, by considering the functions $\theta_i=\sigma_{i+1}-\eta_1$ for all $i=1,\cdots,n$, and using the estimate
$$\sigma_{n+2-p}\geq \frac{n+1-p}{p} \theta_p+ \eta_1\geq  (n+1-p) -\eta_1\frac{(n+1-2p)}{p} \geq n +1-p,$$
we finish the proof with the help of the fact that $\eta_1\leq g(S(\xi),\xi)\leq 0.$
\hfill$\square$\\

\noindent Using this last result, we can prove Corollary \ref{rigidity1} and the next one as in \cite{EHGN1}.
\begin{cor} \label{rigidity2}
Let $N$ be an $(n+2)$-dimensional compact manifold with non-negative curvature operator. Assume that  $M=\mathbb{S}^1 \times \mathbb{S}^n$ with $n\geq 3,$ the sectional curvature $K^N$ of $N$
vanishes on $M,$ the mean curvature $H >0$ and $(n+1-p)\sup_M g(S(\xi),\xi) + 4c^2 [\frac n2]\leq 0$. Then, the manifold $N$ is isometric to $\mathbb{S}^1\times
B'$.
\end{cor}


Also, an analogue result holds: 

\begin{cor}\label{equalitycarapo}
Let $N$ be an $(n+2)$-dimensional compact manifold with non-negative
curvature operator. Assume that the boundary $M$ carries a minimal
Riemannian flow such that $(n+1-p)\sup_M g(S(\xi),\xi) +4 c^2 [\frac n2] \geq 0$ and also admits a
basic special Killing $(n-p)$-form for some $2 \leq p \leq \frac n2$. If the inequality
$\sigma_{p+1}(M)\geq p + \mathop\mathrm{sup}\limits_M g(S(\xi),\xi) + \frac{4}{n+1-p}c^2 [\frac n2]$ holds, the manifold $N$ is isometric to
$\lquot{\mathbb{R} \times B'}{\Gamma}$.
\end{cor}

{\bf Proof.} We follow the same proof as in Corollary \ref{equalitycara}. We just remark that
\begin{eqnarray*}
\sigma_{n+2-p}&\geq&  \frac {n+1-p}{p} (p+\sup_M g(S(\xi),\xi)+\frac{4c^2 [\frac n2]}{n+1-p})+\eta_1(\frac{2p - n-1}{p})\\
&\geq&  (n+1-p)+ \sup_M g(S(\xi),\xi)+\frac{4c^2[\frac n2]}{p}\geq n+1-p. 
\end{eqnarray*}
This finishes the proof of the corollary.
\hfill$\square$\\

\noindent The proof of Corollary \ref{rigidity23} is similar to the one of Corollary  \ref{rigidity1}.

\section{Appendix}
In this section, we will state some technical formulas that we use in our computations. We will omit the proofs of these formulas and will leave them to the reader. \\\\
\noindent For any $X\in \Gamma(TM)$, we have
\begin{equation}\label{eq: appe1}
(\nabla^M_X S)^{[p]}=\nabla^M_XS^{[p]}-S^{[p]}(\nabla^M_X).
\end{equation}
Also, for any $p$-form $\varphi$,
\begin{equation}\label{eq: appe2}
S^{[p+1]}(X\wedge \varphi)=S(X)\wedge \varphi+X\wedge S^{[p]}(\varphi).
\end{equation}
For any orthonormal frame $\{f_i\}_{i=1,\cdots,n+1}$ of $TM$, we have
\begin{equation}\label{eq:appe4}
\langle S^{[p]}\varphi,\varphi\rangle=\sum_{i=1}^{n+1}\langle S(f_i)\lrcorner\varphi,f_i\lrcorner\varphi\rangle.
\end{equation}
Next, we define the interior product of an $s$-form with a $p$-form $\varphi$, the $(p-s)$-form as follows:
$$((X_1\wedge\cdots \wedge X_{s})\lrcorner\varphi)(Y_1,\cdots,Y_{p-s})=\varphi(X_s,\cdots,X_1,Y_1,\cdots,Y_{p-s}).$$
Therefore, the scalar product with any differential $(p-s)$-form satisfies
\begin{equation}\label{eq:appe3}
\langle (X_1\wedge\cdots \wedge X_{s})\lrcorner\varphi,\psi\rangle=(-1)^{\frac{s(s-1)}{2}}\langle \varphi,X_1\wedge\cdots\wedge X_s\wedge\psi\rangle.
\end{equation}
Finally, the action of $S^{[p-s]}$ on $(X_1\wedge\cdots \wedge X_{s})\lrcorner\varphi$ gives the rule
\begin{equation}\label{eq:appe6}
S^{[p-s]}((X_1\wedge\cdots\wedge X_s)\lrcorner \varphi)=(X_1\wedge\cdots\wedge X_s)\lrcorner S^{[p]}(\varphi)-(S^{[s]}(X_1\wedge\cdots\wedge X_s))\lrcorner\varphi.\\
\end{equation}
In particular, this gives for $s=1$
\begin{equation}\label{eq:appe5}
S^{[p-1]}(X\lrcorner \varphi)=X\lrcorner S^{[p]}(\varphi)-S(X)\lrcorner\varphi.
\end{equation}

{\noindent \bf Acknowledgment:} The first two named authors were supported by a fund from the Lebanese University. The second named author would like to thank the Alexander von Humboldt Foundation for its support.


\begin{thebibliography}{99}
\bibitem{BS} M. Belishev and V. Sharafutdinov, \emph{Dirichlet to Neumann operator on differential forms}, Bull. Sci. Math. {\bf 132} (2008), 128-145.
\bibitem{C} Y. Carri\`ere, \emph{Flots riemanniens, structure transverse des feuilletages}, Toulouse, Ast\'erique {\bf 116} (1984), 31-52.
\bibitem{EHGN} F. El Chami, G. Habib, N. Ginoux and R. Nakad, \emph{Rigidity results for spin manifolds with foliated boundary}, to appear in J. Geom.
\bibitem{EHGN1} F. El Chami, G. Habib, O. Makhoul and R. Nakad, \emph{Eigenvalue Estimate for the basic Laplacian on manifolds with foliated boundary}, arXiv:1512.04683, december 2015.


\bibitem{ElG1} A. El Kacimi and B. Gmira, \emph{Stabilit\'e du caract\`ere k\"ahl\'erien transverse}, Israel J. Math. {\bf 101} (1997), 323-347.
\bibitem{ElG2} A. El Kacimi, \emph{Op\'erateurs transversalement \'elliptiques sur un feuilletage riemannien et applications}, Compositio Mathematica {\bf 73} (1990), 57-106.
\bibitem{GM} S. Gallot and D. Meyer, \emph{Sur la premi\`ere valeur propre du $p$-spectre pour les varie\'et\'es \`a op\'erateur de courbure positif}, C.R. Acad. Sci Paris {\bf 320} (1995), 1331-1335.
\bibitem{Li} M. Li, \emph{A sharp comparison theorem for compact manifolds with mean convex boundary}, arxiv:1204.1695v2.
\bibitem{O} B. O'Neill, \emph{The fundamental equations of a submersion}, Mich. Math. J. {\bf 13} (1966), 459-469.
\bibitem{RS} S. Raulot and A. Savo, \emph{A Reilly formula and eigenvalue estimates for differential forms}, J. Geom. Anal. {\bf 3} (2011), 620-640.
\bibitem{RP} K. Richardson and E. Park, \emph{the basic Laplacian of a Riemannian foliation}, Amer. J. Math. {\bf 118} (1996), 1249-1275.
\bibitem{S} G. Schwarz, \emph{Hodge decomposition-A method for solving boundary value problems}, Lecture notes in Mathematics, Springer, 1995.
\bibitem{R} B. Reinhart, \emph{Foliated manifolds with bundle-like metrics}, Ann. Math. {\bf 69} (1959), 119-132.
\bibitem{T} Ph. Tondeur, \emph{Foliations on Riemannian manifolds}, Springer, New York (1959).
\end{thebibliography}
\end{document}